\documentclass[11pt, twoside]{article}
\usepackage{amssymb, amsmath,amsthm}
\usepackage{latexsym}
\usepackage{float,latexsym,shortvrb}
\usepackage[dvips]{graphics, epsfig, color}

\newcommand{\ignore}[1]{}
\ignore{The text written between the parentheses here is ignored by latex}

\newcommand{\R}{\mbox{\rm I}\!\mbox{\rm R}}

\newcommand{\M}[1]{\mathrm{Met}_{#1}}
\newcommand{\m}[1]{\mathrm{met}_{#1}}
\newcommand{\Cu}[1]{\mathrm{Cut}_{#1}}

\newcommand{\cu}[1]{\mathrm{cut}_{#1}}
\newcommand{\Sy}[1]{Sym({#1})}

\def\ds{\delta(S)}

\newcommand{\tot}{910\:209}
\newcommand{\totsubsept}{147\:805}
\newcommand{\totsubhuit}{202\:573}
\newcommand{\totsubNsept}{477}
\newcommand{\totsubNhuit}{389}

\def\@begintheorem#1#2{\par\bgroup{\sc #1\ #2. }\it\ignorespaces}
\def\@opargbegintheorem#1#2#3{\par\bgroup{\sc #1\ #2\ (#3). }\it\ignorespaces}
\def\@endtheorem{\egroup}

\newtheorem{theorem}{Theorem}[section]
\newtheorem{corollary}[theorem]{Corollary}
\newtheorem{lemma}[theorem]{Lemma}
\newtheorem{remark}[theorem]{Remark}
\newtheorem{proposition}[theorem]{Proposition}
\newtheorem{definition}[theorem]{Definition}
\newtheorem{example}[theorem]{Example}
\newtheorem{question}[theorem]{Question}
\newtheorem{problem}[theorem]{Problem}
\newtheorem{conjecture}[theorem]{Conjecture}

\newcommand{\bd}[1]{\begin{definition}\rm\label{#1}}
\newcommand{\bt}[1]{\begin{theorem}\label{#1}}
\newcommand{\bc}[1]{\begin{corollary}\label{#1}}
\newcommand{\bcj}[1]{\begin{conjecture}\label{#1}}
\newcommand{\bl}[1]{\begin{lemma}\label{#1}}
\newcommand{\bp}[1]{\begin{proposition}\label{#1}}
\newcommand{\be}[1]{\begin{example}\rm\label{#1}}
\newcommand{\bq}[1]{\begin{question}\rm\label{#1}}
\newcommand{\bprob}[1]{\begin{problem}\rm\label{#1}}
\newcommand{\beq}[1]{\begin{eqnarray}\label{#1}}
\newcommand{\br}[1]{\begin{remark}\rm\label{#1}}

\newcommand{\el}{\end{lemma}}
\newcommand{\ep}{\end{proposition}}
\newcommand{\ee}{\end{example}}
\newcommand{\eq}{\end{question}}
\newcommand{\eprob}{\end{problem}}
\newcommand{\ecj}{\end{conjecture}}
\newcommand{\eeq}{\end{eqnarray}}
\newcommand{\ed}{\end{definition}}
\newcommand{\et}{\end{theorem}}
\newcommand{\ec}{\end{corollary}}
\newcommand{\er}{\end{remark}}

\begin{document}

 \pagestyle{myheadings}
 \markboth{A counterexample to a conjecture of Laurent and Poljak}
 {Antoine Deza and Gabriel Indik}


\title{\bf A counterexample to a conjecture\\ of Laurent and Poljak}

\author{Antoine  Deza
\and
Gabriel Indik}
\date{November 28, 2005}
\maketitle
\begin{abstract}
The metric polytope $\m{n}$ is the
polyhedron associated with all semimetrics on $n$ nodes and defined by the 
triangle inequalities $x_{ij}-x_{ik}-x_{jk}\leq 0$
and  $x_{ij}+x_{ik}+x_{jk}\leq 2$ for all triples $i,j,k$ of $\{1,\dots,n\}$. 
In 1992   Monique Laurent and Svatopluk Poljak conjectured that every
fractional vertex of the metric polytope is adjacent to some integral vertex.
The conjecture holds for  $n\leq 8$ and, in particular, for the 1 550 825 600 vertices of $\m{8}$.
While the overwhelming majority of the known vertices of $\m{9}$ satisfy the Laurent-Poljak
conjecture, we exhibit a fractional vertex not adjacent to any integral vertex. 
\end{abstract}

\section{Introduction and Notation}
The ${n\choose 2}$-dimensional {\it cut cone} $\Cu{n}$ is usually
introduced as the conic hull of the incidence vectors of all the cuts
of the complete graph on $n$ nodes.  More precisely, given a subset
$S$ of $V_n=\{1,\dots,n\}$, the {\it cut} determined by $S$ consists
of the pairs $(i,j)$ of elements of $V_n$ such that exactly one of
$i$, $j$ is in $S$. By $\ds$ we denote both the cut and its incidence
vector in $\R^{n \choose 2}$, i.e., $\ds_{ij}=1$ if exactly one of
$i$, $j$ is in $S$ and $0$ otherwise for $1\leq i<j\leq n$. 
We use the term cut for both the cut itself and its
incidence vector, so $\ds_{ij}$ are  coordinates of a
point in $\R^{n \choose 2}$. 

The cut cone $\Cu{n}$ is the conic hull of
all $2^{n-1}-1$ nonzero cuts, and the {\it cut polytope} $\cu{n}$ is
the convex hull of all $2^{n-1}$ cuts.  The cut cone and a
relaxation, the {\it metric cone} $\M{n}$, can also be defined in
terms of finite metric spaces in the following way.  For all triples
$\{i,j,k\}\subset V_n$, we consider the following
inequalities.
\begin{eqnarray}
x_{ij}-x_{ik}-x_{jk} &\leq & 0 ,\label{eq:1}\\
x_{ij}+x_{ik}+x_{jk} &\leq & 2 .\label{eq:2}
\end{eqnarray}
(\ref{eq:1}) specify the $3{n \choose 3}$ facets
of the cone $\M{n}$ of {\it semimetrics} on 
$V_n$; that is, of functions 
$x:V_n\times V_n\rightarrow \R_+$
satisfying $x_{ij}=x_{ji}$, $x_{ii}=0$,  and the {\it triangle} 
inequalities (\ref{eq:1}).
While $x$ is a {\it metric} only when $x_{ij}>0$ for all $i\neq j$, we
will follow the usual convention and call $\M{n}$ the {\it metric cone}.

It is well-known that $\Cu{n}$
is the conic hull of all, up to a constant multiple,
$\{0,1\}$-valued extreme rays of $\M{n}$. 
The cuts  satisfy the {\it perimeter} 
inequalities (\ref{eq:2}) which can also be obtained from 
(\ref{eq:1}) by the {\it switching} operation, see Section~\ref{check}.
Bounding $\M{n}$ by the ${n\choose 3}$
facets induced by (\ref{eq:2}), we obtain a natural 
relaxation of $\cu{n}$, the {\it metric polytope}
$\m{n}$, so that
$\cu{n}$ is the convex hull of all $\{0,1\}$-valued vertices of
$\m{n}$.

One of the motivations for the study of these polyhedra comes from
their applications in combinatorial optimization, the most important
being the {\sc maxcut} and multicommodity flow problems.
We refer to Deza and Laurent~\cite{DeLa97} and to Poljak and Tuza~\cite{PT95}
for a detailed study of those
polyhedra and their applications in combinatorial optimization.

\section{A counterexample to the Laurent-Poljak conjecture}
Laurent and Poljak~\cite{LP92} conjectured that every
fractional vertex of the metric polytope $\m{n}$ is adjacent to some integral vertex,
i.e., to a cut. 
Since we have $\m{3}=\cu{3}$ and $\m{4}=\cu{4}$, the conjecture is obviously true
for the $4$ vertices of $\m{3}$ and for the $8$ vertices of $\m{4}$. The conjecture 
holds for the $32$ vertices of $\m{5}$ and the $544$ vertices of $\m{6}$ as well as
for several classes of vertices of $\m{n}$, see~\cite{L96}. The conjecture
was further substantiated by the computation of  $\m{7}$ and $\m{8}$.
The $275\:840$  vertices of $\m{7}$ and the $1\:550\:825\:600$ vertices of $\m{8}$
are adjacent a cut, see~\cite{DDF96,DFMV03,DFPS01}.
 
 While the overwhelming majority of the known vertices of $\m{9}$ satisfy the Laurent-Poljak
conjecture we exhibit a fractional vertex not adjacent to any integral vertex.

\begin{proposition}\label{counter}
The neighbors of the  
fractional vertex $\frac{1}{9}(2$, $2$, $3$, $3$, $4$, $4$, $5$, $5$, $4$, $3$, $5$, $6$, $6$, $3$, $3$, $5$, $5$, $2$, $4$, $3$, $5$, $6$, $3$, $3$, $6$, $6$, $5$, $3$, $2$, $6$, $6$, $3$, $3$, $5$, $3$, $4)$
of the metric polytope 
$\m{9}$
are all fractional.
\end{proposition}
The vertex given in Proposition~\ref{counter}, as well as a few other vertices not adjacent to any cut, were found
 by an extensive computer search of the vertices of the $36$-dimensional metric polytope $\m{9}$, see~\cite{D05}.
Note that while finding a vertex providing a counterexample to the Laurent-Poljak conjecture
is computationally challenging, 
to verify that a given vertex is indeed not adjacent to a cut is
easy if the vertex is quasi-simple, i.e., if the incidence of the given
vertex is equal to the dimension plus one.
For example,  one can easily check, see Section~\ref{check},  that the vertex 
given in Proposition~\ref{counter}
satisfies with equalities $37$ of the  $336$ inequalities defining $\m{9}$ and 
is adjacent to $37$ vertices which are all fractional.

\section{Related questions}
\subsection{The diameter of the metric polytope}
Since any pair of cuts forms an edge of $\m{n}$, the Laurent-Poljak conjecture would imply
that the diameter $\delta(\m{n})$ of the metric polytope satisfies $\delta(\m{n})\leq 3$. 
We recall that the diameter  of a polytope $P$ is the smallest number   $\delta(P)$
such that any two vertices of $P$ can be connected by a path with at most $\delta(P)$ edges.
We have
$\delta(\m{3})=\delta(\m{4})=1$,  $\delta(\m{5})=\delta(\m{6})=2$ and $\delta(\m{7})=\delta(\m{8})=3$.
While the diameter of the restriction of $\m{9}$ to its known vertices appears to be less than $3$, 
it is not clear that the diameter of $\m{n}$ is bounded by a constant.

\subsection{The no-cut set conjecture}
\begin{conjecture}\label{connected}{\rm \cite{DFPS01}}
For $n\geq 6$, the restriction 
of the metric polytope $\m{n}$ to its fractional vertices is connected.
\end{conjecture}

Conjecture~\ref{connected} can be seen as complementary to the Laurent-Poljak conjecture
both graphically and computationally:
For any pair of vertices, while Laurent-Poljak conjecture  implies 
that there is a path made of cuts joining them, i.e., the cut vertices form a {\it dominating set}, Conjecture~\ref{connected} means that there is a path made of non-cut
vertices joining them, i.e., the cut vertices do not form a a {\it cut-set}.
On the other hand, while 
Laurent-Poljak conjecture means that the enumeration of the extreme rays of the metric cone
$\M{n}$ is enough to obtain the vertices of the metric polytope $\m{n}$, Conjecture~\ref{connected}
means that we can obtain the vertices of $\m{n}$ without enumerating the extreme rays of $\M{n}$.

\section{Counterexample generation and verification}\label{check}
One important feature of the metric and cut polyhedra is their very
large symmetry group. We recall that the symmetry group $Is(P)$ of a
polyhedron $P$ is the group of isometries preserving $P$ and that an
isometry is a linear transformation preserving the Euclidean distance.
For $n\geq 5$, the symmetry groups of the polytopes $\m{n}$ and 
$\cu{n}$ are isomorphic and  induced by permutations on
$V_n$ and {\it switching reflections by a cut}, see~\cite{DGL91},
and the symmetry groups of the cones $\M{n}$ and $\Cu{n}$ are isomorphic 
to $\Sy{n}$, see~\cite{DGP03}. Given a
cut $\ds$, the switching reflection $r_{\ds}$ is defined by
$y=r_{\ds}(x)$ where $y_{ij}=1-x_{ij}$ if $(i,j)\in \ds$ and
$y_{ij}=x_{ij}$ otherwise.

\subsection{Counterexample generation}
The vertices of $\m{n}$ are partitioned into orbits under the action of the symmetry group $Is(\m{n})$.
Using a parallel implementation of an orbitwise enumeration algorithm,  $\tot$ orbits of vertices
of  $\m{9}$ were computed on a parallel cluster. 
Among these $\tot$ orbits, $\totsubsept$ are made of vertices belonging to exactly
$37$ inequalities. These vertices are quasi-simple as the dimension of $\m{9}$ is $36$. 
Out of these $\totsubsept$ orbits, $\totsubNsept $ are made of vertices providing a 
counterexample to the 
Laurent-Poljak dominating set conjecture. In addition, out of $\totsubhuit$ 
orbits made of vertices belonging to exactly $38$ inequalities, $\totsubNhuit$  
provide  counterexamples to the Laurent-Poljak conjecture.

\subsection{Counterexample verification}
For a quasi-simple vertex, one can easily verify that all the adjacent vertices are fractional
by performing $3$ elementary computations which we illustrate using the vertex 
given in Proposition~\ref{counter}.
\begin{itemize}
\item[$(i)$]
 Check which of the $336$ inequalities of $\m{9}$  are satisfied with equality by the vertex.
For the vertex given in Proposition~\ref{counter}, we obtain $37$ inequalities, see Section~\ref{Ilist}.
\item[$(ii)$]
Compute the pointed cone formed by the inequalities of $\m{9}$ satisfied with equality by the vertex.
For the vertex given in Proposition~\ref{counter}, we obtain a  quasi-simplicial cone  with $37$ extreme rays.
\item[$(iii)$]
For each extreme ray, perform a ray shooting test from the vertex until piercing one of the facets
 of $\m{9}$ not containing the vertex. For the vertex given in Proposition~\ref{counter},
  we obtain the $37$ fractional vertices given in Section~\ref{Alist}.
\end{itemize}
Note that, while  the computation $(ii)$ can be extremely expensive for a highly degenerate vertex in high dimension, it can be done efficiently if the vertex is quasi-simple. It takes
less than a second of CPU time for the vertex 
given in Proposition~\ref{counter}  using enumeration packages such as   
 {\em lrs}~\cite{LRS}
 or 
 {\em cdd}~\cite{CDD}.
 Computations $(i)$ and $(iii)$ are straightforward and  take less than a second of CPU time.

\subsection{Given counterexample incidence and adjacency lists}\label{lists}
\subsubsection{Given counterexample incidence list}\label{Ilist}
The vertex 
given in Proposition~\ref{counter} satisfies with equalities the following $37$ inequalities
of $\m{9}$:
$\Delta_{6,7,{\bar 9}}$,
$\Delta_{5,{\bar 8},9}$,
$\Delta_{5,{\bar 7},9}$,
$\Delta_{{\bar 5},7,8}$,
$\Delta_{5,6,{\bar 8}}$,
$\Delta_{4,{\bar 7},9}$,
$\Delta_{4,{\bar 6},9}$,
$\Delta_{4,{\bar 6},8}$,
$\Delta_{{\bar 4},6,7}$,
$\Delta_{4,5,9}$,
$\Delta_{4,5,{\bar 7}}$,
$\Delta_{3,{\bar 6},9}$,
$\Delta_{{\bar 3},6,7}$,
$\Delta_{3,5,{\bar 8}}$,
$\Delta_{3,4,{\bar 6}}$,
$\Delta_{2,7,{\bar 9}}$,
$\Delta_{2,6,{\bar 9}}$,
$\Delta_{2,6,{\bar 8}}$,
$\Delta_{2,6,7}$,
$\Delta_{2,5,{\bar 8}}$,
$\Delta_{{\bar 2},4,9}$,
$\Delta_{{\bar 2},4,8}$,
$\Delta_{2,{\bar 4},7}$,
$\Delta_{2,{\bar 4},6}$,
$\Delta_{2,{\bar 3},6}$,
$\Delta_{1,{\bar 5},8}$,
$\Delta_{{\bar 1},4,5}$,
$\Delta_{1,{\bar 3},8}$,
$\Delta_{1,{\bar 3},6}$,
$\Delta_{{\bar 1},3,5}$,
$\Delta_{{\bar 1},3,4}$,
$\Delta_{1,{\bar 2},9}$,
$\Delta_{1,{\bar 2},8}$,
$\Delta_{{\bar 1},2,7}$,
$\Delta_{{\bar 1},2,6}$,
$\Delta_{{\bar 1},2,5}$,
$\Delta_{{\bar 1},2,3}$
where the triangle inequality (\ref{eq:1}) and the perimeter inequality (\ref{eq:2}) are respectively  denoted by
$\Delta_{i,j,{\bar k}}$ 
and 
$\Delta_{i,j,k}$.

\subsubsection{Given counterexample adjacency list}\label{Alist}
The vertex
given in Proposition~\ref{counter} is adjacent to the following $37$ fractional vertices
of $\m{9}$:\\\\
\vspace{1mm} 
$\frac{1}{3}(0, 1, 1, 1, 2, 2, 1, 1, 1, 1, 1, 2, 2, 1, 1, 2, 0, 1, 1, 0, 2, 2, 1, 1, 2, 2, 1, 1, 0, 2, 2, 1, 1, 1, 1, 2)$\\
\vspace{1mm} 
$\frac{1}{3}(0, 1, 1, 1, 2, 2, 1, 1, 1, 1, 1, 2, 2, 1, 1, 2, 2, 1, 1, 2, 2, 2, 1, 1, 2, 2, 1, 1, 0, 2, 2, 1, 1, 1, 1, 2)$\\
\vspace{1mm} 
$\frac{1}{3}(1, 0, 2, 0, 1, 1, 0, 2, 1, 1, 1, 2, 2, 1, 1, 2, 0, 1, 1, 0, 2, 2, 1, 1, 2, 2, 1, 1, 0, 2, 2, 1, 1, 1, 1, 2)$\\
\vspace{1mm} 
$\frac{1}{3}(1, 1, 0, 1, 1, 1, 2, 2, 2, 1, 2, 2, 2, 1, 1, 1, 2, 0, 2, 1, 1, 1, 1, 1, 2, 2, 2, 0, 1, 1, 2, 1, 1, 1, 1, 0)$\\
\vspace{1mm}
$\frac{1}{3}(1, 1, 0, 1, 1, 1, 2, 2, 2, 1, 2, 2, 2, 1, 1, 1, 2, 0, 2, 1, 1, 1, 1, 1, 2, 2, 2, 2, 1, 3, 2, 1, 1, 3, 1, 2)$\\
\vspace{1mm} 
$\frac{1}{3}(1, 1, 1, 1, 1, 1, 2, 2, 1, 1, 2, 2, 2, 1, 1, 2, 2, 1, 1, 1, 2, 2, 1, 1, 2, 2, 2, 1, 1, 2, 2, 1, 1, 2, 1, 1)$\\
\vspace{1mm}
$\frac{1}{3}(1, 1, 1, 1, 1, 1, 2, 2, 2, 1, 2, 2, 2, 1, 1, 1, 2, 0, 2, 1, 1, 2, 1, 1, 2, 2, 2, 1, 1, 2, 2, 1, 1, 2, 1, 1)$\\
\vspace{1mm} 
$\frac{1}{3}(1, 1, 1, 1, 1, 1, 2, 2, 2, 1, 2, 2, 2, 1, 1, 2, 2, 1, 1, 1, 2, 2, 1, 1, 2, 2, 2, 1, 1, 2, 2, 1, 1, 2, 1, 1)$\\
\vspace{1mm}
$\frac{1}{3}(1, 1, 1, 1, 2, 1, 2, 2, 1, 1, 2, 2, 2, 1, 1, 2, 2, 1, 1, 1, 2, 2, 1, 1, 2, 2, 2, 1, 1, 2, 2, 1, 1, 2, 1, 1)$\\
\vspace{1mm} 
$\frac{1}{3}(1, 1, 1, 1, 2, 1, 2, 2, 2, 1, 2, 2, 2, 1, 1, 2, 2, 1, 1, 1, 2, 2, 1, 1, 2, 2, 2, 1, 1, 2, 2, 1, 1, 2, 1, 1)$\\
\vspace{1mm}
$\frac{1}{4}(1, 1, 1, 1, 1, 2, 2, 2, 2, 2, 2, 2, 3, 1, 1, 2, 2, 0, 3, 1, 1, 2, 2, 1, 3, 3, 2, 1, 1, 3, 3, 1, 1, 2, 2, 2)$\\
\vspace{1mm}
$\frac{1}{4}(1, 1, 1, 1, 2, 2, 2, 2, 2, 2, 2, 3, 3, 1, 1, 2, 2, 1, 1, 1, 3, 2, 1, 1, 3, 3, 3, 1, 1, 3, 2, 2, 2, 2, 2, 2)$\\
\vspace{1mm}
$\frac{1}{6}(1, 1, 2, 2, 3, 3, 3, 3, 2, 2, 3, 4, 4, 2, 2, 3, 3, 2, 2, 2, 4, 4, 2, 2, 4, 4, 3, 2, 1, 4, 4, 2, 2, 3, 2, 3)$\\
\vspace{1mm}
$\frac{1}{6}(1, 1, 3, 1, 3, 3, 2, 3, 2, 2, 2, 4, 4, 1, 2, 4, 2, 2, 2, 1, 4, 4, 2, 2, 3, 4, 4, 2, 1, 4, 4, 3, 2, 3, 2, 3)$\\
\vspace{1mm}
$\frac{1}{6}(1, 1, 3, 2, 3, 3, 3, 3, 2, 2, 3, 4, 4, 2, 2, 4, 3, 2, 2, 2, 4, 4, 2, 2, 4, 4, 3, 2, 1, 4, 4, 2, 2, 3, 2, 3)$\\
\vspace{1mm}
$\frac{1}{6}(1, 1, 3, 2, 3, 3, 3, 3, 2, 2, 3, 4, 4, 2, 2, 4, 3, 2, 2, 2, 4, 5, 2, 2, 4, 4, 3, 3, 1, 3, 4, 2, 2, 4, 2, 2)$\\
\vspace{1mm}
$\frac{1}{6}(2, 1, 2, 2, 2, 2, 4, 4, 3, 2, 4, 4, 4, 2, 2, 3, 3, 1, 3, 3, 3, 4, 2, 2, 4, 4, 4, 2, 2, 4, 4, 2, 2, 4, 2, 2)$\\
\vspace{1mm}
$\frac{1}{7}(1, 1, 3, 2, 3, 3, 3, 3, 2, 2, 3, 4, 4, 2, 2, 4, 3, 2, 2, 2, 4, 5, 2, 2, 4, 4, 3, 3, 1, 5, 4, 2, 2, 4, 2, 4)$\\
\vspace{1mm}
$\frac{1}{7}(1, 1, 3, 2, 3, 4, 3, 3, 2, 2, 3, 4, 5, 2, 2, 4, 3, 2, 3, 2, 4, 5, 2, 3, 4, 4, 3, 2, 1, 5, 5, 2, 2, 3, 3, 4)$\\
\vspace{1mm}
$\frac{1}{7}(1, 1, 3, 2, 3, 4, 3, 4, 2, 2, 3, 4, 5, 2, 3, 4, 3, 2, 3, 2, 5, 5, 2, 3, 4, 5, 3, 2, 1, 4, 5, 2, 3, 3, 2, 3)$\\
\vspace{1mm}
$\frac{1}{7}(2, 2, 2, 3, 3, 2, 5, 4, 4, 2, 5, 5, 4, 3, 2, 4, 5, 1, 4, 3, 4, 5, 3, 2, 5, 4, 4, 3, 2, 5, 5, 2, 3, 5, 2, 3)$\\
\vspace{1mm}
$\frac{1}{7}(2, 2, 2, 3, 3, 3, 5, 5, 4, 2, 5, 5, 5, 3, 3, 4, 5, 1, 3, 3, 3, 5, 3, 3, 5, 5, 4, 2, 2, 4, 4, 2, 2, 4, 2, 2)$\\
\vspace{1mm}
$\frac{1}{9}(1, 2, 4, 2, 5, 5, 4, 4, 3, 3, 3, 6, 6, 3, 3, 6, 4, 3, 3, 2, 6, 6, 3, 3, 6, 6, 5, 3, 2, 6, 6, 3, 3, 5, 3, 4)$\\
\vspace{1mm}
$\frac{1}{9}(2, 2, 3, 3, 4, 4, 5, 3, 4, 3, 5, 6, 6, 3, 3, 5, 5, 2, 4, 3, 5, 6, 3, 3, 6, 6, 5, 3, 2, 6, 6, 3, 3, 5, 3, 4)$\\
\vspace{1mm}
$\frac{1}{9}(2, 2, 3, 3, 4, 4, 5, 5, 4, 3, 5, 6, 6, 3, 3, 5, 5, 2, 4, 3, 3, 6, 3, 3, 6, 6, 5, 3, 2, 6, 6, 3, 3, 5, 3, 4)$\\
\vspace{1mm}
$\frac{1}{9}(2, 2, 3, 3, 4, 4, 5, 5, 4, 3, 5, 6, 6, 3, 3, 5, 5, 2, 4, 3, 5, 6, 3, 3, 6, 6, 5, 3, 2, 6, 6, 3, 3, 5, 3, 6)$\\
\vspace{1mm}
$\frac{1}{9}(2, 2, 3, 3, 4, 4, 5, 5, 4, 3, 5, 6, 6, 3, 3, 5, 5, 2, 4, 3, 5, 6, 3, 3, 6, 6, 3, 3, 2, 6, 6, 3, 3, 5, 3, 4)$\\
\vspace{1mm}
$\frac{1}{9}(2, 2, 3, 3, 4, 4, 5, 5, 4, 3, 5, 6, 6, 3, 3, 5, 5, 2, 4, 3, 5, 6, 3, 3, 6, 6, 5, 3, 2, 6, 6, 3, 3, 3, 3, 4)$\\
\vspace{1mm}
$\frac{1}{9}(2, 2, 3, 3, 4, 4, 5, 5, 4, 3, 5, 6, 6, 3, 3, 5, 5, 2, 6, 3, 5, 6, 3, 3, 6, 6, 5, 3, 2, 6, 6, 3, 3, 5, 3, 4)$\\
\vspace{1mm}
$\frac{1}{9}(2, 2, 3, 3, 4, 6, 5, 5, 4, 3, 5, 6, 6, 3, 3, 5, 5, 2, 4, 3, 5, 6, 3, 3, 6, 6, 5, 3, 2, 6, 6, 3, 3, 5, 3, 4)$\\
\vspace{1mm}
$\frac{1}{9}(3, 2, 2, 4, 3, 3, 6, 6, 5, 3, 5, 6, 6, 3, 3, 4, 6, 1, 5, 4, 4, 6, 3, 3, 6, 6, 5, 3, 2, 6, 6, 3, 3, 5, 3, 4)$\\
\vspace{1mm}
$\frac{1}{10}(2, 2, 4, 3, 5, 4, 5, 5, 4, 4, 5, 7, 6, 3, 3, 6, 5, 3, 4, 3, 7, 7, 3, 4, 7, 7, 6, 3, 2, 6, 7, 4, 4, 5, 3, 4)$\\
\vspace{1mm}
$\frac{1}{10}(2, 2, 4, 3, 5, 5, 5, 6, 4, 4, 5, 7, 7, 3, 4, 6, 5, 3, 3, 3, 6, 7, 3, 3, 7, 6, 6, 4, 2, 7, 6, 4, 3, 6, 3, 5)$\\
\vspace{1mm}
$\frac{1}{10}(3, 3, 2, 4, 4, 3, 7, 7, 6, 3, 7, 7, 6, 4, 4, 5, 7, 1, 6, 4, 4, 6, 4, 3, 7, 7, 6, 3, 3, 7, 7, 3, 3, 6, 4, 4)$\\
\vspace{1mm}
$\frac{1}{12}(3, 3, 3, 5, 5, 5, 7, 7, 6, 4, 8, 8, 8, 4, 4, 6, 8, 2, 6, 4, 6, 8, 4, 4, 8, 8, 8, 4, 4, 8, 8, 4, 4, 8, 4, 4)$\\
\vspace{1mm}
$\frac{1}{12}(3, 3, 3, 5, 5, 5, 8, 7, 6, 4, 8, 8, 8, 5, 4, 6, 8, 2, 6, 5, 6, 8, 4, 4, 7, 8, 6, 4, 3, 8, 8, 3, 4, 7, 4, 5)$\\
\vspace{1mm}
$\frac{1}{12}(3, 3, 3, 5, 5, 5, 8, 7, 6, 4, 8, 8, 8, 5, 4, 6, 8, 2, 6, 5, 6, 8, 4, 4, 9, 8, 8, 4, 3, 8, 8, 5, 4, 7, 4, 5)$\\\\

\noindent
{\bf acknowledgements}
Research supported  by the 
Natural Sciences and Engineering Research Council of Canada
under the Canada Research Chair and the Discovery Grant  programs.
Thanks to the Shared Hierarchical Academic Research Computing Network
(SHARCNET) for a generous allocation of CPU time.\\

\noindent
{\small Antoine Deza, Gabriel Indik}\\
Advanced Optimization Laboratory,
Department of Computing and Software,\\
McMaster University, Hamilton, Ontario, Canada. \\
{\em Email}: deza, indikg{\small @}mcmaster.ca.

\end{document}